\magnification=1200
\input amstex
\input amssym.def
\documentstyle{amsppt}
\parindent=0pt
\font\fabs=cmr8 scaled\magstep0

\define\Z{{\Cal Z}}

\define\z+{{\Bbb Z}_+}

\define\A_n{{\Cal A}_q(n)}

\define\U_n{{\Cal U}_q(n)}

\define\ep{\varepsilon}

\define\wi{\widetilde}

\define\ov{\overline}
\define\hb{\hfill\break}

\redefine\a{\alpha}
\redefine\C{{\Bbb C}}
\redefine\b{\beta}

\redefine\d{\delta}
\define\dpr{{^{\prime\hskip-0.5pt\prime}}}
\define\pr{{^{\prime}}}
\define\ot{\otimes}
\define\sq{$\hfill \square$}
\headline={\ifodd\pageno\rightheadline\else\leftheadline\fi}
\def\rightheadline{\tenrm\hfil A Non-Commutative Discrete Hypergroup
Associated with $q$-Disk Polynomials \hfil\folio}
\def\leftheadline{\tenrm\folio\hfil Paul G.A. Floris \hfil}
\voffset=2\baselineskip

\topmatter
\title A Non-Commutative Discrete Hypergroup Associated With $q$-Disk
Polynomials\endtitle
\author Paul G.A. Floris\endauthor
\affil Leiden University \\ Department of Mathematics and Computer Science \\
P.O. Box 9512 \\ 2300 RA Leiden\\ The Netherlands\\ e-mail: floris\@wi.leidenuniv.nl\\ fax: +31 - 71 - 27 69 85 \endaffil
\endtopmatter

{\fabs {{\sl Keywords and phrases} : $q$-disk polynomials, linearization
coefficients, DJS-hypergroup.}}
\bigskip
{\bf 0. Introduction}.
\medskip
In connection with a given system of orthogonal polynomials $\{p_n\}$
it is of great interest to know if there exist positive measures
$\mu_{x,y}(z)$ and non-negative numbers $c_{k,l}(m)$ such that
$$
p_n(x) p_n(y)=\int p_n(z) d\mu_{x,y}(z) \tag 0.1
$$
and
$$
p_k(x) p_l(x)=\sum_m c_{k,l}(m) p_m(x). \tag 0.2
$$
The first formula, called product formula, gives rise to a positive 
convolution structure and
the second formula, called linearization formula, to a dual positive 
convolution structure associated
with these orthogonal polynomials.
It is a quite classical result that such positivity results hold for
Gegenbauer polynomials $P_n^{(\alpha,\alpha)}$ ($\alpha\ge-1/2$).
Around 1970 these results were also proved for more general Jacobi
polynomials $P_n^{(\alpha,\beta)}$ with $(\alpha,\beta)$ in a set
containing $\{(\alpha,\beta)\in{\Bbb R}^2\mid \alpha\ge\beta\ge-1/2\}$
(see [G1] and [G2]).
In case such polynomials have an interpretation as spherical functions
on a compact symmetric space $G/K$, these positivity results and
associated convolution structure follow immediately from analysis on
the space of $K$-biinvariant functions on $G$. In connection with
this, see also the survey paper by Gasper [G3]. In the seventies the
essential aspects of such zonal analysis on groups were abstracted into
the concept of a DJS-hypergroup by work of Dunkl, Jewett and Spector.
This made it possible to conclude that the positivity results in the
two product formulas for Jacobi polynomials with quite general
parameters $\alpha, \beta$ give rise to two associated hypergroups,
one discrete and one continuous, and dual to each other.
For an elaborate exposition of the theory of hypergroups we
refer the reader to the book by Bloom and Heyer [BH].\hb
A formula like $(0.1)$ can often be given in an explicit way. Then it may
also be possible to extend it, using Carlson's theorem, by analytic
continuation from some discrete set of parameter values to a more general
set. However, in $(0.2)$ the coefficients are often not explicitly known.
Thus the analytic continuation method will not work, and one must look for
an alternative way of proving positivity for the more general parameter
set. One way of achieving this uses an addition formula
for the polynomials $p_n$ for the general set of parameters, obtained
from the discrete case by a continuation argument (see the method
described in [Koo1]).\hb
In the theory of quantum groups it is a natural question to ask whether one 
can obtain results similar to the ones we have in the classical situation.
In his paper [Koo4]
Koornwinder showed that it is possible to associate a discrete hypergroup
with the 'double coset space' of a
Gel'fand pair of compact quantum groups, although the construction is
somewhat more involved than the classical one. But as in the classical 
situation the basic ingredient is positivity of
linearization coefficients for the related spherical functions. 
It is perhaps good to note that the 
hypergroups arising in this way need not be commutative.
\smallskip
The aim of this paper is to give an example of a non-commutative
discrete hypergroup associated with $q$-disk polynomials. These are
polynomials $R_{l,m}^{(\a)}$ in two non-commuting variables which are
expressed through little $q$-Jacobi polynomials and that appear, 
for the value $\a=n-2$, as zonal spherical 
functions on a quantum analogue of the homogeneous space $U(n)/U(n-1)$.
This fact was first proved in [NYM] (see also [Fl]). In a previous paper
[Fl] we proved an addition formula for these $q$-disk polynomials. It is
this addition formula that will allow us to prove positivity of linearization
coefficients in a manner similar to [Koo1], and to construct from it a
DJS-hypergroup following [Koo4].
\smallskip
The  paper is organized as follows. In section one we briefly
recall the definition of $q$-disk polynomials and some of their properties.
Furthermore we will state the addition formula which they satisfy. Section
two merely deals with the proof of positivity, or rather non-negativity,
of the linearization coefficients. The proof resembles the way of reasoning
in [Koo1]. Finally, in section three we explicitly
construct the non-commutative discrete hypergroup related to the $q$-disk
polynomials.
\medskip
We end by fixing the notation and recalling some well-known
facts. In all that follows we will keep $0<q<1$ fixed.\hb
Recall the definition of the little $q$-Jacobi polynomials:
$$
p_m(x;a,b;q) = {_2}\varphi_1 \bigl[ {{q^{-m}\ ,\ abq^{m+1}}\atop {aq}}; 
q,qx\bigr] = \sum_{k=0}^m {{(q^{-m};q)_k (abq^{m+1};q)_k}\over {(aq;q)_k 
(q;q)_k}} (qx)^k.
$$
If $0<aq<1$ and $bq<1$ they satisfy the orthogonality
$$
\sum_{k=0}^\infty {{(bq;q)_k}\over{(q;q)_k}}(aq)^k \big( p_lp_m\bigr)
(q^k;a,b;q) =\d_{l,m} {{(q,bq;q)_l (aq)^l (1-abq) (abq^2;q)_\infty}\over
{(aq,abq;q)_l (1-abq^{2l+1}) (aq;q)_\infty}}.
$$
Here 
$$
\eqalign{
(a;q)_{\infty} &= \prod_{j=0}^{\infty} (1-q^ja)\cr
(a;q)_{\b} &= {{(a;q)_{\infty}}\over {(aq^{\b};q)_{\infty}}}\cr}
$$
are $q$-shifted factorials $(\b\in\C)$.
In particular, if we let
$$
P_m ^{(\a,\b)} (x;q) = p_m(x;q^{\a},q^{\b};q)\qquad\quad (\a,\b>-1)
$$
then the orthogonality reads
$$
\eqalign{
\int_0^1 P_l^{(\a,\b)}(x;q) &P_m^{(\a,\b)}(x;q) x^{\a} 
{{(q x;q)_{\infty}}\over {(q^{\b+1}x;q)_{\infty}}} d_{q} x =\cr
&\d_{lm}{{(1-q)q^{m(\a+1)}}\over {1-q^{\a+\b+2m+1}}} {{(q;q)_m (q;q)_{\b+m}}\over 
{(q^{\a+1};q)_m (q^{\a+1};q)_{\b+m}}}.\cr}
$$
Here we used Jackson's $q$-integral
$$
\int_0^c f(x)d_q x = c(1-q) \sum_{k=0}^{\infty} f(cq^k) q^k.
$$
Finally, write $\z+$ for the non-negative integers:
$$
\z+= \{0,1,2,\hdots\}.
$$
{\bf Acknowledgement }: The author thanks Professor Tom H. Koornwinder for
his useful suggestions.
\bigskip
{\bf 1. $q$-Disk polynomials and their addition formula.}
\medskip
Suppose we are given a complex unital $\ast$-algebra $\Cal Z$ generated by 
the elements $z$ and $z^{\ast}$, subject to the relation
$$
z^{\ast} z= q^2 zz^{\ast} + 1-q^2\tag 1.1
$$
and with $\ast$-structure $(z)^{\ast}=z^{\ast}$. 
It is not hard to show that $\Z$ has as a linear basis the set
$\{ z^k (z^{\ast})^l\ :\ k,l\in\z+\}$.\hb
On this
algebra we define the {\it $q$-disk polynomials} $R_{l,m}^{(\a)} 
(z,z^{\ast};q)$ for $\a>-1$ and $l,m\in \z+$ as follows:
$$
R_{l,m}^{(\a)}(z,z^{\ast};q) =\left\{\aligned z^{l-m} P_m^{(\a,l-m)}
(1-zz^{\ast};q) \qquad (l\geq m)\\
P_l^{(\a,m-l)}(1-zz^{\ast};q) (z^{\ast})^{m-l}\qquad (l\leq m).
\endaligned\right.\tag 1.2
$$
Note that
$$
R_{l,m}^{(\a)}(z,z^{\ast};q)^{\ast} = R_{m,l}^{(\a)}(z,z^{\ast};q).\tag 1.3
$$
It is easily seen that one has $R_{l,m}^{(\a)}(z,z^{\ast};q)=
\sum_{i=0}^{l\wedge m} c_i z^{l-i}(z^{\ast})^{m-i}$ with $c_0\neq 0$.
The orthogonality of these polynomials can be expressed through $q$-integrals:
$$
\eqalign{
{1\over {2\pi}}\int_0^1 \int_0^{2\pi}
R_{l,m}^{(\a)}(e^{i\theta}z,e^{-i\theta}z^{\ast};q^2)^{\ast} 
R_{l\pr,m\pr}^{(\a)}&(e^{i\theta}z,e^{-i\theta}z^{\ast};q^2)
 d \theta\cr 
(1&-zz^{\ast})^{\a} d_{q^2}(1-zz^{\ast}) \cr
= \d_{ll\pr}\d_{mm\pr} &{{1-q^2}\over
{1-q^{2(\a+1)}}} c_{l,m}^{(\a)}\cr}\tag 1.4
$$
where
$$
c_{l,m}^{(\a)} = {{(1-q^{2(\a+1)})q^{2m(\a+1)}}\over {1-q^{2(\a+l+m+1)}}} {{(q^2;q^2)_l 
(q^2;q^2)_{m}}\over {(q^{2(\a+1)};q^2)_l (q^{2(\a+1)};q^2)_{m}}}.\tag 1.5
$$
Note that $(1.4)$ is well-defined, since after integrating with respect
to $\theta$ one obtains a polynomial which is invariant under the transformation
$z\mapsto e^{i\theta}z,\ z^{\ast}\mapsto e^{-i\theta}z^{\ast}$, and hence is
a polynomial in the single variable $1-zz^{\ast}$.\hb
The same orthogonality can
be achieved using the linear functional $h_{(\a)} : \Z\to \C$ $(\a>-1)$
defined
as
$$
h_{(\a)}(z^k (z^{\ast})^l) = \d_{kl} q^{2k(\a +1)} {{(q^2;q^2)_k}\over
{(q^{2(\a +2)};q^2)_k}}
$$
and satisfying $h_{(\a)}(p^{\ast}) = \ov{h_{(\a)}(p)}$ for all $p\in \Z$.
Then 
$$
h_{(\a)}\bigl(R_{l,m}^{(\a)}(z,z^{\ast};q^2)^{\ast} 
R_{l\pr,m\pr}^{(\a)}(z,z^{\ast};q^2)\bigr) =\d_{ll\pr}\d_{mm\pr}
c_{l,m}^{(\a)}.\tag 1.6
$$
\hb
It follows that the $R_{l,m}^{(\a)}(z,z^{\ast};q)\
(l,m\in\z+)$ form an orthogonal basis for $\Z$ with respect to the inner
product defined by $h_{(\a)}$. For more details we refer the reader to [Fl].
\medskip
In [Fl, Thm. 3.5.8] we proved the
following {\it addition formula} for these $q$-disk polynomials:

\proclaim\nofrills {\bf Theorem 1} : Suppose we are given the abstract complex 
$\ast$-algebras $\Cal X$ and $\Cal Y$ with generators $X_1,X_2,X_1^{\ast},
X_2^{\ast}$ and $Y_1,Y_2,Y_1^{\ast},
Y_2^{\ast}$ respectively, relations 
$$
\eqalign{
X_1X_2 = &qX_2X_1\cr
X_1^{\ast}X_2 = &qX_2X_1^{\ast}\cr
X_2^{\ast}X_2 = &q^2X_2X_2^{\ast} + (1-q^2)\cr
X_1^{\ast}X_1 = &q^2 X_1X_1^{\ast}  
+(1-q^2) (1-X_2X_2^{\ast})\cr
Y_1Y_2 = &qY_2Y_1\cr
Y_1^{\ast}Y_2 = &qY_2Y_1^{\ast}\cr
Y_1Y_1^{\ast} = &Y_1^{\ast}Y_1\cr
1 =&Y_1Y_1^{\ast} + Y_2Y_2^{\ast}
=q^2Y_1^{\ast}Y_1 + Y_2^{\ast}Y_2.\cr}
\tag 1.7
$$
and $\ast$-structures
$$
\alignat 3
&(X_1)^{\ast} = X_1^{\ast}  & &\qquad & 
&(Y_1)^{\ast}= Y_1^{\ast}\tag 1.8\\
&(X_2)^{\ast}= X_2^{\ast} & &\qquad & 
&(Y_2)^{\ast}= Y_2^{\ast}.
\endalignat
$$
Then, for arbitrary $\a>0$ and arbitrary $l,m\in\z+$ we have the following 
addition formula for the $q$-disk polynomials:
$$
\eqalign{
R_{l,m}^{(\a)}(-q X_1\ot &Y_1^{\ast}+X_2\ot Y_2, -qX_1^{\ast}\ot Y_1
+X_2^{\ast}\ot Y_2^{\ast};q^2) =\cr
\sum_{r=0}^l\sum_{s=0}^m c_{l,m;r,s}^{(\a)} &R_{l-r,m-s}^{(\a+r+s)}
(X_2,X_2^{\ast}; q^2)
R_{r,s}^{(\a-1)}(X_1,X_1^{\ast}, 1-X_2X_2^{\ast}; q^2)\cr
&\ot
(-q)^{r-s} R_{l-r,m-s}^{(\a+r+s)} (Y_2,Y_2^{\ast}; q^2) Y_1^s 
(Y_1^{\ast})^r. \cr}
\tag 1.9
$$
\endproclaim
Here
$$
c_{l,m;r,s}^{(\a)} = {{1-q^{2(\a+r+s+1)}}\over  {1-q^{2(\a+1)}}} 
{{c_{l,m}^{(\a)}}\over {c_{l-r,m-s}^{(\a+r+s)} 
c_{r,s}^{(\a-1)}}} 
$$
(cf. (1.5)), and
$$
R_{l,m}^{(\a)}(A,B,C;q) =\left\{\aligned C^m A^{l-m} P_m^{(\a,l-m)}({{C-AB}
\over C};q)
\qquad (l\geq m)\\
C^l P_l^{(\a,m-l)}({{C-AB}\over C};q) B^{m-l}\qquad (l\leq m).
\endaligned\right.\tag 1.10
$$
With this notation we have that $R_{l,m}^{(\a)}(A,B,1;q) =
R_{l,m}^{(\a)}(A,B;q)$.
\bigskip
Define the following map on $\Cal X$:
$$
\wi{h}_{(\a)} (p(X_2, X_2^{\ast}) X_1^k (X_1^{\ast})^l) =
\d_{kl}\ p(X_2, X_2^{\ast}) (1-X_2X_2^{\ast})^k 
q^{2k\a} {{(q^2;q^2)_k}\over {(q^{2(\a +1)};q^2)_k}}
$$
where $p(X_2, X_2^{\ast})$ is any (ordered) polynomial in $X_2, X_2^{\ast}$,
and $k,l\in\z+$.
\smallskip
\proclaim\nofrills {\bf Lemma 2}: Let $p_1(X_2, X_2^{\ast})$ and 
$p_2(X_2, X_2^{\ast})$ be arbitrary ordered polynomials in 
$X_2, X_2^{\ast}$ and let $p_3(X_1, X_1^{\ast}, 1- X_2 X_2^{\ast})$
be any ordered polynomial in $X_1, X_1^{\ast}, 1- X_2 X_2^{\ast}$,
homogeneous of degree $k$, where we put $deg (X_1) = deg (X_1^{\ast})
= {1\over 2}$ and $deg(1- X_2 X_2^{\ast})=1$. Then
$$
\eqalign{
&\wi{h}_{(\a)} (p_1(X_2, X_2^{\ast}) p_3(X_1, X_1^{\ast}, 1- X_2 X_2^{\ast})
p_2(X_2, X_2^{\ast})) =\cr
&h_{(\a -1)}(p_3(X_1, X_1^{\ast}, 1))\ p_1(X_2, X_2^{\ast})
(1-X_2X_2^{\ast})^k p_2(X_2, X_2^{\ast}).\cr}
$$
\endproclaim
{\sl Proof }: In view of $(1.7)$ we can write
$$
 p_3(X_1, X_1^{\ast}, 1- X_2 X_2^{\ast}) =\sum_{i=0}^k\sum_{j=0}^{k-i}
c_{ij} (1-X_2X_2^{\ast})^i X_1^j (X_1^{\ast})^{-j +2(k-i)}.
$$
The assertion now easily follows when one uses the first two lines of $(1.7)$ 
and the relations 
$X_2(1-X_2X_2^{\ast}) =
q^{-2} (1-X_2X_2^{\ast})X_2$ and $X_2^{\ast} (1-X_2X_2^{\ast}) =
q^2 (1-X_2X_2^{\ast})X_2^{\ast}$, which are immediate from 
$(1.7)$. \sq
\medskip
Combining Lemma 2 with $(1.3)$, $(1.6)$ and $(1.5)$ we obtain as a 
consequence
\proclaim\nofrills {\bf Corollary 3}: We have
$$
\eqalign{
\wi{h}_{(\a)} (&R_{l-i,m-j}^{(\a +i+j)}(X_2,X_2^{\ast};q^2)
R_{i,j}^{(\a -1)}(X_1,X_1^{\ast}, 1-X_2X_2^{\ast};q^2)\times\cr
&R_{p,r}^{(\a -1)}(X_1,X_1^{\ast}, 1-X_2X_2^{\ast};q^2)^{\ast}
R_{l\pr-p,m\pr-r}^{(\a +p+r)}(X_2,X_2^{\ast};q^2)^{\ast})=\cr
\d_{ip}\d_{jr} &c_{j,i}^{(\a -1)} 
R_{l-i,m-j}^{(\a +i+j)}(X_2,X_2^{\ast};q^2) (1-X_2X_2^{\ast})^{i+j}
R_{l\pr-p,m\pr-r}^{(\a +p+r)}(X_2,X_2^{\ast};q^2)^{\ast}.\cr}
$$
\endproclaim

\bigskip

{\bf 2. Positivity of linearization coefficients.}
\medskip
Since the $R_{l,m}^{(\a)}(z,z^{\ast};q^2)$ $(l,m\in\z+)$ form a basis for 
$\Z$, we have the following expansion in $\Z$:
$$
\eqalign{
&R_{l,m}^{(\a)}(z,z^{\ast};q^2) R_{l\pr,m\pr}^{(\a)}(z,z^{\ast};q^2)^{\ast} 
=\cr
&\sum_{l\dpr,m\dpr\in \z+} a(l,m;l\pr,m\pr;l\dpr,m\dpr)
R_{l\dpr,m\dpr}^{(\a)}(z,z^{\ast};q^2).\cr}\tag 2.1
$$
Here only finitely many of the coefficients
$a(l,m;l\pr,m\pr;l\dpr,m\dpr)$ are non-zero if we fix $l,m,l\pr,m\pr$.
In fact the sum ranges over those values of $l\dpr, m\dpr$ such that
$l\dpr-m\dpr = l-m -(l\pr-m\pr)$.
The $a(l,m;l\pr,m\pr;l\dpr,m\dpr)$ are called {\it linearization
coefficients}.
\smallskip
\proclaim\nofrills {\bf Theorem 4}: For all $\a>0$ and all possible choices
of $(l,m), (l\pr,m\pr), (l\dpr, m\dpr)\in \z+^2$ the linearization coefficients
are non-negative:
$$
a(l,m;l\pr,m\pr;l\dpr,m\dpr)\geq 0.
$$
\endproclaim
{\sl Proof }: First note that both the pair $X_2,X_2^{\ast}$ and the pair
$Y_2,Y_2^{\ast}$ satisfy $(1.1)$.
Let $\Omega = -qX_1\ot Y_1^{\ast} + X_2\ot Y_2$. It is 
straightforward to verify that $\Omega^{\ast}\Omega = q^2 \Omega\Omega^{\ast}
+ 1-q^2$. This means that we have an identity similar to (2.1) but with
$z, z^{\ast}$ replaced by $\Omega, \Omega^{\ast}$:
$$
\eqalign{
&R_{l,m}^{(\a)}(\Omega, \Omega^{\ast};q^2) R_{l\pr,m\pr}^{(\a)}
(\Omega, \Omega^{\ast};q^2)^{\ast} =\cr
&\sum_{l\dpr,m\dpr\in \z+} a(l,m;l\pr,m\pr;l\dpr,m\dpr)
R_{l\dpr,m\dpr}^{(\a)}(\Omega, \Omega^{\ast};q^2).\cr}\tag 2.2
$$
Now substitute $(1.9)$ into the right-hand side as well as for both factors
in the left-hand side of $(2.2)$ and apply $\wi{h}_{(\a)}\ot id$ to this.
By Lemma 2 and $(1.6)$ we obtain for the right-hand side
$$
\eqalign{
\sum_{l\dpr,m\dpr} &\sum_{r=0}^{l\dpr}\sum_{s=0}^{m\dpr} 
a(l,m;l\pr,m\pr;l\dpr,m\dpr) c_{l\dpr,m\dpr;r,s}^{(\a)} \times\cr
&\wi{h}_{(\a)} \biggl( R_{l\dpr-r,m\dpr-s}^{(\a+r+s)}(X_2,X_2^{\ast}; q^2)
R_{r,s}^{(\a-1)}(X_1,X_1^{\ast}, 1-X_2X_2^{\ast}; q^2)\biggr) \ot\cr
&(-q)^{r-s} R_{l\dpr-r,m\dpr-s}^{(\a+r+s)} (Y_2,Y_2^{\ast}; q^2) Y_1^s 
(Y_1^{\ast})^r \cr
\sum_{l\dpr,m\dpr}
&a(l,m;l\pr,m\pr;l\dpr,m\dpr)
c_{l\dpr,m\dpr;0,0}^{(\a)} c_{0,0}^{(\a -1)}\times\cr
&R_{l\dpr,m\dpr}^{(\a)}
(X_2, X_2^{\ast},1; q^2) \ot
R_{l\dpr,m\dpr}^{(\a)} (Y_2,Y_2^{\ast},1; q^2) =\cr
\sum_{l\dpr,m\dpr} &a(l,m;l\pr,m\pr;l\dpr,m\dpr)\times\cr
&R_{l\dpr,m\dpr}^{(\a)}
(X_2, X_2^{\ast},1; q^2) \ot
R_{l\dpr,m\dpr}^{(\a)} (Y_2,Y_2^{\ast},1; q^2).\cr}
$$
On the left we get
$$
\eqalign{
\sum_{i=0}^l&\sum_{j=0}^m \sum_{p=0}^{l\pr}\sum_{r=0}^{m\pr}
c_{l,m;i,j}^{(\a)} c_{l\pr,m\pr;p,r}^{(\a)} (\wi{h}_{(\a)}\ot id)\cr
\biggl( &R_{l-i,m-j}^{(\a+i+j)}(X_2, X_2^{\ast};q^2) R_{i,j}^{(\a-1)} 
(X_1,X_1^{\ast},1-X_2X_2^{\ast};q^2)\times\cr
&R_{p,r}^{(\a-1)} (X_1,X_1^{\ast},1-X_2X_2^{\ast};q^2)^{\ast}
R_{l\pr-p,m\pr-r}^{(\a+p+r)}(X_2, X_2^{\ast};q^2)^{\ast}\ot\cr
&(-q)^{i-j} R_{l-i,m-j}^{(\a+i+j)}(Y_2,Y_2^{\ast},1;q^2) (Y_1^{\ast})^i 
Y_1^j\times\cr
&(-q)^{p-r}(Y_1^{\ast})^r Y_1^p R_{l\pr-p,m\pr-r}^{(\a+p+r)}(Y_2,Y_2^{\ast};q^2)^{\ast}
\biggr )=\cr
\sum_{i=0}^{l\wedge l\pr} &\sum_{j=0}^{m\wedge m\pr}  
c_{l,m;i,j}^{(\a)} c_{l\pr,m\pr;i,j}^{(\a)}  
c_{j,i}^{(\a-1)} \times \cr
&R_{l-i,m-j}^{(\a+i+j)}(X_2,X_2^{\ast};q^2)
(1-X_2X_2^{\ast})^{i+j} 
R_{l\pr-i,m\pr-j}^{(\a+i+j)}(X_2,X_2^{\ast};q^2)^{\ast}
\ot\cr
&R_{l-i,m-j}^{(\a+i+j)}(Y_2,Y_2^{\ast};q^2) 
(1-Y_2Y_2^{\ast})^{i+j} 
R_{l\pr-i,m\pr-j}^{(\a+i+j)}(Y_2,Y_2^{\ast};q^2)^{\ast} \cr}
$$
in view of Corollary 3.
So by now we have the identity
$$
\eqalign{
\sum_{i=0}^{l\wedge l\pr} &\sum_{j=0}^{m\wedge m\pr}  
c_{l,m;i,j}^{(\a)} c_{l\pr,m\pr;i,j}^{(\a)}  
c_{j,i}^{(\a-1)} \times \cr
&R_{l-i,m-j}^{(\a+i+j)}(X_2,X_2^{\ast};q^2)
(1-X_2X_2^{\ast})^{i+j}
R_{l\pr-i,m\pr-j}^{(\a+i+j)}(X_2,X_2^{\ast};q^2)^{\ast}
\ot\cr
&R_{l-i,m-j}^{(\a+i+j)}(Y_2,Y_2^{\ast};q^2) 
(1-Y_2Y_2^{\ast})^{i+j} R_{l\pr-i,m\pr-j}^{(\a+i+j)}(Y_2,Y_2^{\ast};q^2)^{\ast}
=\cr
\sum_{l\dpr,m\dpr} &a(l,m;l\pr,m\pr;l\dpr,m\dpr)
R_{l\dpr,m\dpr}^{(\a)}
(X_2, X_2^{\ast}; q^2) \ot
R_{l\dpr,m\dpr}^{(\a)} (Y_2,Y_2^{\ast}; q^2).\cr}\tag 2.3
$$
Write $\sigma$ for the $\ast$-algebra anti-automorphism $\sigma :\Cal Y
\to \Cal Y$ which interchanges $Y_2$ and $Y_2^{\ast}$ and fixes $Y_1$ and
$Y_1^{\ast}$. Note that
$$
\sigma \bigl( R_{r,s}^{(\a)}(Y_2,Y_2^{\ast};q^2) \bigr) = 
R_{r,s}^{(\a)}(Y_2,Y_2^{\ast};q^2)^{\ast}.
$$
Letting $id \ot\sigma$ act on $(2.3)$ yields
$$
\eqalign{
\sum_{l\dpr,m\dpr} a(l,&m;l\pr,m\pr;l\dpr,m\dpr) 
R_{l\dpr,m\dpr}^{(\a)}(X_2, X_2^{\ast}; q^2) \ot
R_{l\dpr,m\dpr}^{(\a)} (Y_2,Y_2^{\ast}; q^2)^{\ast}=\cr
\sum_{i=0}^{l\wedge l\pr} \sum_{j=0}^{m\wedge m\pr}  
&c_{l,m;i,j}^{(\a)} c_{l\pr,m\pr;i,j}^{(\a)}  
c_{j,i}^{(\a-1)} \times \cr
&R_{l-i,m-j}^{(\a+i+j)}(X_2,X_2^{\ast};q^2)
(1-X_2X_2^{\ast})^{i+j}
R_{l\pr-i,m\pr-j}^{(\a+i+j)}(X_2,X_2^{\ast};q^2)^{\ast}
\ot\cr
&R_{l\pr-i,m\pr-j}^{(\a+i+j)}(Y_2,Y_2^{\ast};q^2) 
(1-Y_2Y_2^{\ast})^{i+j} R_{l-i,m-j}^{(\a+i+j)}(Y_2,Y_2^{\ast};q^2)^{\ast}.\cr}
$$
Finally, multiply from the left with $R_{l\dpr,m\dpr}^{(\a)}
(X_2, X_2^{\ast}; q^2)^{\ast} \ot 1$ and from the right with
$1\ot R_{l\dpr,m\dpr}^{(\a)} (Y_2,Y_2^{\ast}; q^2)$
and evaluate $h_{(\a)}\ot h_{(\a)}$ on the result. By virtue of $(1.6)$ we
wind up with:
$$
\eqalign{
a(l,m;l\pr,m\pr;&l\dpr,m\dpr) (c_{l\dpr, m\dpr}^{(\a)})^2
=
\sum_{i=0}^{l\wedge l\pr} \sum_{j=0}^{m\wedge m\pr}  
c_{l,m;i,j}^{(\a)} c_{l\pr,m\pr;i,j}^{(\a)}  
c_{j,i}^{(\a-1)} \times \cr
\biggl| h_{(\a)} \biggl(&R_{l\dpr,m\dpr}^{(\a)}
(X_2, X_2^{\ast}; q^2)^{\ast} R_{l-i,m-j}^{(\a+i+j)}(X_2,X_2^{\ast};q^2)\cr
&(1-X_2X_2^{\ast})^{i+j}
R_{l\pr-i,m\pr-j}^{(\a+i+j)}(X_2,X_2^{\ast};q^2)^{\ast}\biggr)
\biggr| ^2\cr}
$$
because $h_{(\a)}$ satisfies $h_{(\a)}(p^{\ast}) = \ov{h_{(\a)}(p)}$.
Since $0<q<1$, this will imply that $a(l,m;l\pr,m\pr;l\dpr,m\dpr)\geq 0$.\sq
\medskip
{\bf Remark }: Considering the case where $l=m$ we thus obtain non-negativity 
for the linearization coefficients of the little $q$-Jacobi polynomials
$P_m^{(\a,0)}(x)$ $(\a>-1)$.
\bigskip

{\bf 3. A discrete hypergroup structure associated with $q$-disk
polynomials.}
\medskip
In this section we construct a so-called DJS-hypergroup from the
linearization formula $(2.1)$ in the way it was pointed out by Koornwinder 
[Koo4].
First we define the proper setting.\hb
Let $K$ be a locally compact Hausdorff topological space and let $M(K)$
be the space of all complex regular Borel measures on $K$, and $M^1(K)$
the subset of all probability measures. For $x\in K$ we denote by $\d_x$
the corresponding point measure: $\d_x (x) =1$ (so $\d_x\in M^1(K)$).
Assume that in addition there exist
\parindent=20pt
\item{(a)} {\it convolution} : a continuous map $K\times K\to M^1(K),\
(x,y)\to \d_x\star\d_y$ in the weak topology with respect to $C_c(K)$
\item{(b)} {\it involution} : an involutive homeomorphism $K\to K,\ x\to 
{\ov x}$
\item{(c)} {\it unit element}: a distinguished element $e\in K$.
\parindent=0pt
\smallskip
Upon identifying $x$ with $\d_x$, the map in (a) extends uniquely to a
continuous bilinear map $M(K)\times M(K)\to M(K),\ (\mu,\nu)\to \mu\star \nu$,
The involution of (b) induces an involution $\mu\to\mu^{\ast}$ on $M(K)$
as follows: $\mu^{\ast}(E) = \ov {\mu(\overline E)}$ ($E\subset K$ a
Borel subset).
\medskip
{\bf Definition }: The quadruple $(K,\star, {}^-,e)$ is called a
{\it DJS-hypergroup} if for all $x,y,z\in K$ the following conditions are met:
\parindent=20pt
\item{(1)} $\d_x\star (\d_y\star \d_z) = (\d_x\star\d_y)\star\d_z$ 
\item{(2)} supp$(\d_x\star\d_y)$ is compact 
\item{(3)} $(\d_x\star\d_y)^{\ast} = \d_{\ov y}\star\d_{\ov x}$
\item{(4)} $\d_e\star\d_x = \d_x =\d_x\star \d_e$
\item{(5)} $e\in \text{supp} (\d_{\ov x}\star \d_y)$ if and only if $x=y$
\item{(6)} the map of $K\times K$ to the space of nonvoid compact subsets of $K$
given by $(x,y)\to \text{supp} (\d_x\star\d_y)$ is continuous. Here the target
space has the topology as defined in [Je, \S 2.5].
\parindent=0pt
\smallskip
The hypergroup is called {\it commutative} if $\d_x\star \d_y=\d_y\star \d_x$
for all $x,y\in K$, otherwise it is called {\it non-commutative}.
\smallskip
\proclaim\nofrills {\bf Theorem 5}: Put $K= \z+^2$, endowed with the discrete
topology. For $(l,m), (l\pr,m\pr)$ and $(l\dpr,m\dpr)\in K$ define
$$
(\d_{(l,m)}\star\d_{(l\pr,m\pr)})((l\dpr,m\dpr)) =
a(l,m;l\pr,m\pr;l\dpr,m\dpr)
$$
with $a(l,m;l\pr,m\pr;l\dpr,m\dpr)$ as in (2.1). As an involution on $K$
take $(l,m)^- =(m,l)$. Furthermore write $e =(0,0)$.\hb
Then the quadruple
$(K,\star, {^-},e)$ forms a non-commutative discrete DJS-hypergroup.
\endproclaim
{\sl Proof }: We have to verify (a), (1)-(6). Let us abbreviate 
$R_{l,m}^{(\a)}(z,z^{\ast};q^2)$ by $R_{l,m}^{(\a)}$. \hb
(a): The convolution is continuous since we have given $K$ the discrete
topology. Define the multiplicative linear functional $\ep : \Z\to \C$
by $\ep(z) =1= \ep(z^{\ast})$ (so in fact $\ep$ is identically $1$ on
$\Z$). If we now apply $\ep$ to both side of $(2.1)$ we get
$$
1 = \sum_{l\dpr,m\dpr} a(l,m;l\pr,m\pr;l\dpr,m\dpr) 
$$
hence $\d_{(l,m)}\star\d_{(l\pr,m\pr)}\in M^1(K)$.\hb
(1): From $R_{l,m}^{(\a)}(R_{l\pr,m\pr}^{(\a)} R_{l\dpr,m\dpr}^{(\a)})=
(R_{l,m}^{(\a)}R_{l\pr,m\pr}^{(\a)}) R_{l\dpr,m\dpr}^{(\a)}$ it follows
that
$$
\sum_{r,s} a(l,m;l\pr,m\pr;r,s) a(r,s;l\dpr,m\dpr;u,v)=
\sum_{r,s} a(l\pr,m\pr;l\dpr,m\dpr;r,s) a(l,m;r,s;u,v),
$$
whence $\d_{(l,m)}\star (\d_{(l\pr,m\pr)}\star \d_{(l\dpr,m\dpr)}) = 
(\d_{(l,m)}\star \d_{(l\pr,m\pr)})\star \d_{(l\dpr,m\dpr)}$.\hb
(2): Since only finitely many of the elements $a(l,m;l\pr,m\pr;l\dpr,m\dpr)$
are non-zero when $(l,m), (l\pr,m\pr)$ are fixed, the support of 
$\d_{(l,m)}\star\d_{(l\pr,m\pr)}$ is compact.\hb
(3): Using $(1.3)$ we see that $(R_{l,m}^{(\a)} R_{l\pr,m\pr}^{(\a)})^{\ast}=
R_{m\pr,l\pr}^{(\a)} R_{m,l}^{(\a)}$ which gives
$$
\sum_{l\dpr,m\dpr} a(l,m;l\pr,m\pr;m\dpr,l\dpr)R_{l\dpr,m\dpr}
^{(\a)}=
\sum_{l\dpr,m\dpr} a(m\pr,l\pr;m,l;l\dpr,m\dpr) R_{l\dpr,m\dpr}
^{(\a)}.
$$
From this we obtain that 
$$
\eqalign{
(\d_{(l,m)}\star\d_{(l\pr,m\pr)})^{\ast} &= 
\d_{(m\pr,l\pr)}\star\d_{(m,l)}\cr
&=\d_{(l\pr,m\pr)^-}\star\d_{(l,m)^-}\cr}
$$
(4): Since $R_{0,0}^{(\a)}=1$ we get that $a(l,m;0,0;m\dpr,l\dpr)=
\d_{l\dpr 0}\d_{m\dpr 0}$.\hb
(5): Note that
$$
\eqalign{
\d_{ll\pr}\d_{mm\pr} c_{l,m}^{(\a)} &= h_n ({R_{l,m}^{(\a)}}^{\ast}
R_{l\pr,m\pr}^{(\a)})\cr
&= h_n (R_{m,l}^{(\a)} R_{l\pr,m\pr}^{(\a)})\cr
&= \sum_{l\dpr,m\dpr} a(m,l;l\pr,m\pr;m\dpr,l\dpr)
h_n( R_{l\dpr,m\dpr}^{(\a)})\cr
&= a(m,l;l\pr,m\pr;0,0).\cr}
$$
So $e=(0,0)\in \text{supp} (\d_{(l,m)^-}\star\d_{(l\pr,m\pr)})=
\text{supp} (\d_{(m,l)}\star\d_{(l\pr,m\pr)})$ if and only if
one has $a(m,l;l\pr,m\pr;0,0)\neq 0$, which, by the above, 
is true only in case $(l,m)=(l\pr,m\pr)$.\hb
(6): Obvious since $K$ has discrete topology.\sq
\bye